\documentclass{amsart}
\usepackage{latexsym}

\usepackage{amsfonts,amssymb,amsmath}
\usepackage{amsthm}
\usepackage{graphicx,color}
\usepackage{color}
\usepackage{amsmath} 
\usepackage{multirow,bigdelim}
\usepackage{cleveref}
\newtheorem{theorem}{Theorem}	
\crefname{theorem}{Theorem}{Theorems}
\newtheorem{lemma}{Lemma}
\newtheorem{corollary}{Corollary}		
\newtheorem{proposition}{Proposition}	
\crefname{proposition}{Proposition}{Propositions}

\newtheorem{definition}{Definition}

\newtheorem{remark}{Remark}

\crefname{section}{Section}{Sections}
\crefname{theorem}{Theorem}{Theorems}
\crefname{lemma}{Lemma}{Lemmas}
\crefname{corollary}{Corollary}	{Corollaries}			
\crefname{proposition}{Proposition}{Propositions}	
\crefname{claim}{Claim}{Claims}
\crefname{conjecture}{Conjecture}{Conjectures}			
\crefname{definition}{Definition}{Definitions}
\crefname{problem}{Problem}{Problems}
\crefname{example}{Example}{Examples}
\crefname{remark}{Remark}{Remarks}
\crefname{figure}{Figure}{Figures}

\newfont{\bg}{cmr9 scaled\magstep2}

\makeatletter
\newcommand{\tpitchfork}{%
  \vbox{
    \baselineskip\z@skip
    \lineskip-.52ex
    \lineskiplimit\maxdimen
    \m@th
    \ialign{##\crcr\hidewidth\smash{$-$}\hidewidth\crcr$\pitchfork$\crcr}
  }%
}

\newcommand{\R}{\mathbb{R}}

\makeatother
\title[Generic linear
perturbations]
{
Generic linear perturbations
}
\author{Shunsuke Ichiki
}
\thanks{Research Fellow DC1 of Japan Society for the Promotion of Science}
\address{
Graduate School of Environment and Information Sciences,  
Yokohama National University, 
Yokohama 240-8501, Japan}
\email{ichiki-shunsuke-jb@ynu.jp}

\begin{document}
\date{}
\begin{abstract}
In his celebrated paper ``Generic projections'', 
John Mather has shown that 
almost all linear projections from a 
submanifold of a vector space into a subspace 
are transverse with respect to a given modular submanifold. 
In this paper, an improvement of Mather's result is stated. 
Namely, we show that almost all linear perturbations of a smooth mapping 
from a submanifold of $\mathbb{R}^m$ into $\mathbb{R}^\ell$ 
yield a transverse mapping with respect to a given modular submanifold. 
Moreover, applications of this result are given. 
\end{abstract}
\maketitle
\noindent
\section{Introduction}\label{section 1}

Throughout this paper, $\ell$, $m$, $n$ stand for positive integers. 
In this paper, unless otherwise stated, all manifolds and mappings belong to class $C^{\infty}$ 
and all manifolds are without boundary. 

An $n$-dimensional manifold is 
denoted by $N$. 
Let $\pi:\mathbb{R}^m\to \mathbb{R}^\ell$ be a linear mapping. 

In a celebrated paper \cite{GP}, for a given embedding $f:N\to \mathbb{R}^m$, 
a composition $\pi\circ f:N\to \mathbb{R}^\ell$ $(m>\ell)$ 
is investigated, and 
the following assertions (M1)-(M5) 
are obtained for a generic mapping. 
All of (M1)-(M5) follow from the main result 
(\cref{@m} in \cref{section 2}) proved 
by Mather. 
\begin{enumerate}
\item[(M1)] If $(n,\ell)=(n,1)$, 
then a generic function $\pi\circ f:N\to \mathbb{R}$ is a Morse function. 
\item[(M2)] If $(n,\ell)=(2,2)$,
then a generic mapping $\pi\circ f:N\to \mathbb{R}^2$ is an excellent map in the sense defined 
by Whitney in \cite{plane to plane}. 
\item[(M3)] If $(n,\ell)=(2,3)$, 
then the only singularities of the image of a generic mapping 
$\pi\circ f:N\to \mathbb{R}^3$ are 
normal crossings and pinch points.
\item[(M4)] A generic mapping $\pi\circ f:N\to \mathbb{R}^{\ell}$ 
is transverse with respect to the Thom-Boardman 
varieties (for the definition of Thom-Boardman 
varieties, refer to \cite{Arnold}, \cite{Boardman}, \cite{TB}, \cite{Thom}). 
\item[(M5)]If $(n,\ell)$ is in the nice range of dimensions 
(for the definition of nice range of dimensions, refer to \cite{mather6}),  
then a generic mapping $\pi\circ f:N\to \mathbb{R}^{\ell}$ 
is locally infinitesimally stable (for the definition of local infinitesimal stability, 
see \cref{section 2}). 
If moreover, $N$ is compact, 
then a generic mapping $\pi\circ f:N\to \mathbb{R}^{\ell}$ 
is stable (for the definition of stability, \
see \cref{section 2}).
\end{enumerate}
Let $\mathcal{L}(\mathbb{R}^{m},\mathbb{R}^{\ell})$ 
be the space consisting of linear mappings 
of $\mathbb{R}^{m}$ into $\mathbb{R}^{\ell}$. 
For a given embedding $f:N \to \mathbb{R}^m$, 
a property of mappings $\pi\circ f:N\to \mathbb{R}^\ell$ will be said to be 
true for a {\it generic mapping} 
if there exists a subset $\Sigma$ with Lebesgue measure zero of 
$\mathcal{L}(\mathbb{R}^{m},\mathbb{R}^{\ell})$ 
such that for any $\pi\in \mathcal{L}(\mathbb{R}^{m},\mathbb{R}^{\ell})-\Sigma$, 
$\pi\circ f:N\to \mathbb{R}^{\ell}$ has the property. 

The main aim of this paper is to 
prove \cref{@f} in \cref{section 2}, which is an improvement of \cref{@m} 
in \cref{section 2}, 
shown by Mather (\cite{GP}). 

Let $U$ be 
an open subset of $\mathbb{R}^m$ and 
let $F:U\to \mathbb{R}^{\ell}$ be a mapping. 
For any $\pi\in\mathcal{L}(\mathbb{R}^{m},\mathbb{R}^{\ell})$, 
set $F_{\pi}$ as follows: 
\begin{eqnarray*}
F_{\pi}=F+\pi. 
\end{eqnarray*}
Here, the mapping $\pi$ in $F_{\pi}=F+\pi$ is restricted to $U$. 
For a given embedding $f:N\to U$ and 
a given mapping $F:U\to \mathbb{R}^\ell$, 
a property of mappings $F_{\pi}\circ f:N\to \mathbb{R}^\ell$ will be said to be 
true for a {\it generic mapping} 
if there exists a subset $\Sigma$ with Lebesgue measure zero of 
$\mathcal{L}(\mathbb{R}^{m},\mathbb{R}^{\ell})$ 
such that for any $\pi\in \mathcal{L}(\mathbb{R}^{m},\mathbb{R}^{\ell})-\Sigma$, 
$F_\pi \circ f:N\to \mathbb{R}^{\ell}$ has the property. 
For a given embedding $f:N\to U$, from \cref{@f}, 
the following assertions (I1)-(I5) hold.  
\begin{enumerate}
\item[(I1)] If $(n,\ell)=(n,1)$, 
then a generic function $F_{\pi}\circ f:N\to \mathbb{R}$ is a Morse function. 
\item[(I2)] If $(n,\ell)=(2,2)$, 
then a generic mapping $F_{\pi}\circ f:N\to \mathbb{R}^2$ is an excellent map. 
\item[(I3)]If $(n,\ell)=(2,3)$, 
then the only singularities of the image of 
a generic mapping $F_{\pi}\circ f:N\to \mathbb{R}^3$ are 
normal crossings and pinch points.
\item[(I4)] A generic mapping $F_{\pi}\circ f:N\to \mathbb{R}^{\ell}$ 
is transverse with respect to the Thom-Boardman varieties.
\item[(I5)]If $(n,\ell)$ is in the nice range of dimensions, 
then a generic mapping $F_{\pi}\circ f:N\to \mathbb{R}^{\ell}$ 
is locally infinitesimally stable. 
If moreover, $N$ is compact, 
then a generic mapping $F_{\pi}\circ f:N\to \mathbb{R}^{\ell}$ is stable. 
\end{enumerate}
The assertion (M5) (resp., (I5)) above implies assertions 
(M1), (M2), and (M3) (resp., (I1),  (I2)  and (I3)). 
Both assertions (M4) and (M5) (resp., (I4) and (I5)) follow from 
\cref{@m} (resp., \cref{@f}) of \cref{section 2}. 
Moreover, in the special case of $F=0$, $U=\mathbb{R}^{m}$ and $m>\ell$, 
the assertions (I1)-(I5) 
are exactly the same as the assertions (M1)-(M5) 
respectively. 
Note that in the case $m\leq\ell$, 
a generic mapping $\pi \circ f:N \to \mathbb{R}^{\ell}$ is an embedding. 
On the other hand, in the same case, 
a generic mapping $F_\pi \circ f:N \to \mathbb{R}^{\ell}$ is not 
necessarily an embedding. 

The original motivation for this work is to investigate 
the stability of 
quadratic mappings of $\mathbb{R}^m$ into 
$\mathbb{R}^\ell$ of a special type 
called ``generalized distance-squared mappings'' 
(for the precise definition of generalized distance-squared mappings, 
see \cref{section 4}). 
In \cite{G1} (resp., \cite{G2}), 
the generalized distance-squared mappings 
in the case $(m,\ell)=(2,2)$ (resp., $(m,\ell)=(k+1,2k+1)$) 
have been investigated, where $k$ is a positive integer. 
As an application of (I5), 
if $(m,\ell)$ is in the nice range of dimensions, 
then it is shown that 
a generic generalized distance-squared mapping of 
$\mathbb{R}^m$ into $\mathbb{R}^\ell$ 
is locally infinitesimally stable 
(see \cref{@Glis} and \cref{remark G} in \cref{section 4}). 

Notice that for example, the references \cite{brucekirk} and \cite{semi} are also important papers related to generic projections.  
In \cite{brucekirk}, 
an improvement of Mather's result is given 
by replacing a given embedding $f:N\to \mathbb{R}^m$ 
by a given stable mapping $f:N\to \mathbb{R}^m$ (see Theorem 2.2 in 
\cite{brucekirk}). 
On the other hand, in this paper, 
an improvement of Mather's result is given 
by replacing generic projections 
by generic linear perturbations. 

\par 
\bigskip 
In \cref{section 2}, some standard definitions and 
the important notion of ``modular'' submanifold (\cref{modular}) 
defined in \cite{GP} are reviewed, and 
the main theorem (\cref{@f}) in this paper is stated. 
\cref{section 3} is devoted to 
the proof of \cref{@f}. 
In \cref{section 4}, 
the motivation to investigate the stability of 
generalized distance-squared mappings 
is given in detail, and 
as applications of the main theorem, 
results containing \cref{@Glis} are stated. 
\section{Preliminaries and the statement of the main result}\label{section 2}
Let $N$ and $P$ be manifolds and let $J^r(N,P)$ be the space of 
$r$-jets of mappings of $N$ into $P$. 
For a given mapping $g:N\to P$, 
the mapping $j^r g:N\to J^r(N,P)$ 
is defined by $q \mapsto j^r g(q)$. 
Let $C^\infty(N,P)$ be the set of $C^\infty$ mappings of $N$ into $P$, 
and the topology on $C^\infty(N,P)$ is the Whitney $C^\infty$ topology 
(for the definition of Whitney $C^\infty$ topology, 
see for example \cite{GG}). 
Given $g, h\in C^\infty(N,P)$, 
we say that $g$ is {\it $\mathcal{A}$-equivalent} to $h$ 
if there exist diffeomorphisms $\Phi:N\to N$ and 
$\Psi:P\to P$ such that $g=\Psi \circ h \circ \Phi^{-1}$. 
Then, $g$ is said to be {\it stable} 
if the $\mathcal{A}$-equivalence class of $g$ is open in $C^\infty(N,P)$. 

Let $s$ be a positive integer. Define $N^{(s)}$ as follows:
\begin{eqnarray*}
N^{(s)}=\{(q_1,\ldots ,q_s) \mid q_i\not =q_j  (1\leq i<j\leq s) \}. 
\end{eqnarray*}
Let ${}_s J^r(N,P)$ be the space consisting of 
elements $(j^r g(q_1),\ldots ,j^r g(q_s))\in J^r(N,P)^s$ 
satisfying $(q_1,\ldots ,q_s)\in N^{(s)}$. 
Since $N^{(s)}$ is an open submanifold of $N^s$, 
the space ${}_s J^r(N,P)$ is also an open submanifold of $J^r(N,P)^s$.
For a given mapping $g:N\to P$, 
the mapping ${}_s j^rg:N^{(s)}\to {}_s J^r(N,P)$ 
is defined by $(q_1,\ldots ,q_s) \mapsto (j^rg(q_1),\ldots ,j^rg(q_s))$.  

Let $W$ be a submanifold of ${}_s J^r(N,P)$. 
For a given mapping $g:N\to P$, 
we say that ${}_s j^r g:N^{(s)}\to {}_s J^r(N,P)$ is {\it transverse} to $W$ 
if for any $q\in N^{(s)}$, ${}_s j^r g(q)\not\in W$ or  
in the case of ${}_s j^r g(q)\in W$, the following holds: 
\begin{eqnarray*}
d({}_s j^r g)_q(T_qN^{(s)})+T_{{}_s j^r g(q)}W=T_{{}_s j^r g(q)}{}_s J^r(N,P).
\end{eqnarray*}
A mapping $g:N\to P$ will be said to be {\it transverse with respect to} $W$ 
if ${}_s j^r g:N^{(s)}\to {}_s J^r(N,P)$ is transverse to $W$. 

Following Mather (\cite{GP}), we can partition $P^s$ as follows. 
Given any partition $\pi$ of $\{1,\ldots ,s\}$, 
let $P^{\pi}$ denote the set of $s$-tuples 
$(y_1,\ldots ,y_s)\in P^s$ 
such that $y_i=y_j$ 
if and only if 
two positive integers 
$i$ and $j$ are in the same member of the partition $\pi$. 

Let Diff $N$ denote 
the group of diffeomorphisms of $N$. 
There is a natural action of Diff $N$ $\times $ Diff $P$ 
on ${}_s J^r(N,P)$ 
such that for a mapping $g:N\to P$, 
the equality 
$(h,H)\cdot {}_sj^rg(q)={}_sj^r(H\circ g\circ h^{-1})(q')$ 
holds, where $q=(q_1,\ldots ,q_s)$ and 
$q'=(h(q_1),\ldots ,h(q_s))$.  
A subset $W$ of ${}_s J^r(N,P)$ is said to be {\it invariant} 
if it is invariant under this action. 

We recall the following identification $(\ast)$ from \cite{GP}. 
Let $q=(q_1,\ldots ,q_s)\in N^{(s)}$, 
let $g:U\to P$ be a mapping defined 
in a neighborhood $U$ of $\{q_1,\ldots ,q_s\}$ in $N$, and 
let $z={}_s j^r g(q)$, $q'=(g(q_1),\ldots ,g(q_s))$. 
Let ${}_s J^r(N,P)_q$ and ${}_s J^r(N,P)_{q,q'}$ denote 
the fibers of ${}_s J^r(N,P)$ over $q$ and over $(q,q')$ respectively. 
Let $J^r(N)_q$ denote the $\mathbb{R}$-algebra of 
$r$-jets at $q$ of functions on $N$. 
Namely, 
\begin{eqnarray*}
J^r(N)_q={}_s J^r(N,\mathbb{R})_q.
\end{eqnarray*} 
Set $g^*TP=\bigcup_{\widetilde{q}\in U}T_{g(\widetilde{q})}P$, 
where $TP$ is the tangent bundle of $P$.    
Let $J^r(g^*TP)_q$ denote 
the $J^r(N)_q$-module of $r$-jets at $q$ of 
sections of the bundle $g^*TP$. 
Let $\mathfrak{m}_q$ be the ideal in $J^r(N)_q$ 
consisting of jets of functions which vanish at $q$. 
Namely,   
\begin{eqnarray*}
\mathfrak{m}_q=\{{}_s j^r h(q)\in {}_s J^r(N,\mathbb{R})_q \mid 
h(q_1)=\cdots =h(q_s)=0\}. 
\end{eqnarray*}
Let $\mathfrak{m}_q J^r(g^*TP)_q$ be the set consisting of finite sums 
of products of 
an element of $\mathfrak{m}_q$ and an element of $J^r(g^*TP)_q$. 
Namely, we get 
\begin{eqnarray*}
\mathfrak{m}_q J^r(g^*TP)_q
=J^r(g^*TP)_q\cap \{{}_sj^r \xi(q)\in {}_s J^r(N,TP)_q \mid \xi(q_1)=\cdots =\xi(q_s)=0\}.
\end{eqnarray*}
Then, it is seen that the following canonical identification 
of $\mathbb{R}$ vector spaces $(\ast)$ holds.
\begin{align}
T({}_s J^r(N,P)_{q,q'})_z=\mathfrak{m}_q J^r(g^*TP)_q \tag{$\ast $}.
\end{align}

Let $W$ be a non-empty submanifold of 
${}_s J^r(N,P)$. 
Choose $q=(q_1,\ldots ,q_s)\in N^{(s)}$ and 
$g:N\to P$, and 
let $z={}_s j^r g(q)$ and $q'=(g(q_1),\ldots ,g(q_s))$. 
Suppose that the choice is made so that $z\in W$. 
Set $W_{q,q'}=\widetilde{\pi}^{-1}(q,q')$, where 
$\widetilde{\pi}:W\to N^{(s)}\times P^s$ is defined by 
$\widetilde{\pi}({}_s j^r \widetilde{g}(\widetilde{q}))
=(\widetilde{q}, (
\widetilde{g}(\widetilde{q}_1), \ldots ,
\widetilde{g}(\widetilde{q}_s)))$ and 
$\widetilde{q}=(\widetilde{q}_1,\ldots , \widetilde{q}_s)\in N^{(s)}$. 
Suppose that $W_{q,q'}$ is a submanifold of  ${}_s J^r(N,P)$. 
Then, under the identification $(\ast)$, the tangent space 
$T(W_{q,q'})_z$ can be identified with a vector subspace of 
$\mathfrak{m}_q$$J^r(g^*TP)_q$. 
We denote this vector subspace by $E(g,q,W)$. 

\begin{definition}\label{modular}
{\rm 
We say that a submanifold $W$ of  ${}_s J^r(N,P)$ is {\it modular} 
if conditions $(\alpha)$ and $(\beta)$ below are satisfied:
\begin{enumerate}
\item[$(\alpha)$]
The set $W$ is an invariant submanifold of ${}_s J^r(N,P)$, 
and lies over $P^{\pi}$ 
for some partition $\pi$ of $\{1,\ldots ,s\}$.
\item[$(\beta)$] 
For any $q\in N^{(s)}$ and any mapping $g:N \to P$ 
such that ${}_s j^r g(q)\in W$, 
the subspace $E(g,q,W)$ is a $J^r(N)_q$-submodule.
\end{enumerate}
}
\end{definition}

Now, suppose that $P=\mathbb{R}^\ell$. 
The main theorem of \cite{GP} is the following. 
\begin{theorem}[\cite{GP}] \label{@m}
Let $N$ be a manifold of dimension $n$. 
Let $f$ be an embedding 
of $N$ into $\mathbb{R}^m$. 
If $W$ is a modular submanifold of ${}_{s}J^{r}(N,\mathbb{R}^\ell)$ and $m>\ell$, 
then there exists a subset $\Sigma$ with Lebesgue measure zero 
of $\mathcal{L}(\mathbb{R}^{m},\mathbb{R}^{\ell})$ 
such that for any $\pi \in \mathcal{L}(\mathbb{R}^{m},\mathbb{R}^{\ell})-\Sigma$, 
$\pi \circ f:
N\to \mathbb{R}^\ell$ is transverse with respect to $W$.
\end{theorem}

Then, the main theorem in this paper is the following. 
\begin{theorem}\label{@f}
Let $N$ be a manifold of dimension $n$. 
Let $f$ be an embedding 
of $N$ into an open subset $U$ of $\mathbb{R}^m$. 
Let $F:U\to \mathbb{R}^\ell$ be a mapping. 
If $W$ is a modular submanifold of ${}_{s}J^{r}(N,\mathbb{R}^\ell)$, 
then there exists a subset $\Sigma$ with Lebesgue measure zero 
of $\mathcal{L}(\mathbb{R}^{m},\mathbb{R}^{\ell})$ 
such that for any $\pi \in \mathcal{L}(\mathbb{R}^{m},\mathbb{R}^{\ell})-\Sigma$, 
$F_\pi \circ f:
N\to \mathbb{R}^\ell$ is transverse with respect to $W$.
\end{theorem}

It follows that the Thom-Boardman varieties are modular by Mather 
(see \cite{TB} and \cite{GP}). 
Hence, we have the following as a corollary of \cref{@f}. 
\begin{corollary}
\label{F TB}
Let $N$ be a manifold of dimension $n$. 
Let $f$ be an embedding 
of $N$ into an open subset $U$ of $\mathbb{R}^m$. 
Let $F:U\to \mathbb{R}^\ell$ be a mapping. 
Then, there exists a subset $\Sigma$ with Lebesgue measure zero 
of $\mathcal{L}(\mathbb{R}^{m},\mathbb{R}^{\ell})$ 
such that 
for any $\pi \in \mathcal{L}(\mathbb{R}^{m},\mathbb{R}^{\ell})-\Sigma$,  
$F_\pi \circ f:
N\to \mathbb{R}^\ell$ is transverse with respect to 
the Thom-Boardman varieties.
\end{corollary} 

Let $S$ be a finite subset of $N$ and $y$ be a point of $P$. 
Let $g:(N,S)\to (P,y)$ be a map-germ. 
A map-germ $\xi:(N,S)\to (TP, \xi(S))$ such that $\Pi \circ \xi=g$ 
is called a {\it vector field along $g$}, 
where $\Pi : TP\to P$ is the canonical projection. 
Let $\theta(g)_{S}$ be the set consisting of vector fields along $g$. 
Set $\theta(N)_{S}=\theta({\rm id}_N)_{S}$ and $\theta(P)_{y}=\theta({\rm id}_P)_{y}$, 
where ${\rm id}_N:(N,S)\to (N,S)$ and ${\rm id}_P:(P,y)\to (P,y)$ 
are the identify map-germs. 
The mapping $tg:\theta(N)_{S}\to \theta(g)_{S}$ is defined by $tg(\xi)=Tg \circ \xi$, 
where $Tg:TN\to TP$ is the derivative mapping of $g$. 
The mapping $\omega g :\theta(P)_{y}\to \theta(g)_{S}$ is defined by 
$\omega g (\eta) =\eta \circ g$. 
Then, a mapping $g:N\to P$ is said to be {\it locally infinitesimally stable} 
if the following holds 
for every $y\in P$ and every finite subset $S\subset g^{-1}(y)$
(\cite{mather5} and \cite{GP}). 
\begin{eqnarray*}
tg(\theta(N)_{S})+\omega g(\theta(P)_{y})=\theta(g)_{S}.
\end{eqnarray*}

By the same way as in the proof of Theorem 3 of \cite{GP}, 
we have the following as a corollary of \cref{@f}.
\begin{corollary}
\label{F locally infinitesimally stable}
Let $N$ be a manifold of dimension $n$. 
Let $f$ be an embedding 
of $N$ into an open subset $U$ of $\mathbb{R}^m$. 
Let $F:U\to \mathbb{R}^\ell$ be a mapping. 
If a dimension pair $(n,\ell)$ is in the nice dimensions, 
then there exists a subset $\Sigma$ with Lebesgue measure zero 
of $\mathcal{L}(\mathbb{R}^{m},\mathbb{R}^{\ell})$ 
such that 
for any $\pi \in \mathcal{L}(\mathbb{R}^{m},\mathbb{R}^{\ell})-\Sigma$, 
the composition $F_{\pi }\circ f:N\to \mathbb{R}^\ell$ 
is locally infinitesimally stable. 
\end{corollary} 

\begin{remark}\label{remark F}
{\rm 
\begin{enumerate}
\item 
In the special case that  
$F=0$, 
$U=\mathbb{R}^{m}$, and $m>\ell $, \cref{@f} is \cref{@m}. 
\item 
The set $\Sigma$ in Mather's theorem (\cref{@m}) depends only on 
$f:N\to \R^m$ and a modular submanifold $W$ of ${}_{s}J^{r}(N,\mathbb{R}^\ell)$. 
On the other hand,  $\Sigma$ in the main theorem of this paper (\cref{@f}) 
depends on $F:U\to \R^\ell$ too. 
\item
Suppose that the mapping 
$F_{\pi }\circ f:N\to \mathbb{R}^{\ell}$ is proper 
in \cref{F locally infinitesimally stable}. 
Then, the local infinitesimal stability of $F_{\pi }\circ f$ 
implies the stability of it (see \cite{mather5}). 
\item
We explain the advantage that the domain of the mapping $F$ 
is not $\R^m$ but an open set $U$.  
Suppose that $U=\mathbb{R}$. 
Let $F:\mathbb{R}\to \mathbb{R}$ be the mapping defined by $x\mapsto |x|$. 
Since $F$ is not differentiable at $x=0$, we can not apply \cref{@f} to 
the mapping $F:\mathbb{R}\to \mathbb{R}$. 

On the other hand, if $U=\mathbb{R}-\{0\}$, then \cref{@f} 
can be applied to the restriction $F|_{U}$.
\end{enumerate}
}
\end{remark}
\section{Proof of \cref{@f}}\label{section 3}
Let $(\alpha_{ij})_{1\leq i \leq \ell, 1\leq j \leq m}$ be a representing matrix of 
a linear mapping $\pi:\mathbb{R}^m\to \mathbb{R}^\ell$. 
Set $F_{\alpha}=F_{\pi}$, and we have
\begin{eqnarray*}
F_{\alpha}(x)=
\biggl(F_{1}(x)+\sum_{j=1}^{m}\alpha_{1j}x_j, \ldots,
F_{\ell}(x)+\sum_{j=1}^{m}\alpha_{\ell j}x_j\biggr), 
\end{eqnarray*}
where $\alpha=(\alpha_{11}, \ldots ,\alpha_{1m},\ldots ,\alpha_{\ell 1},\ldots ,\alpha_{\ell m})\in (\mathbb{R}^m)^\ell$, $F=(F_1, \ldots ,F_\ell)$ and $x=(x_1,\ldots ,x_m)$. 
For a given embedding $f:N\to U$, 
a mapping $F_{\alpha}\circ f:N\to \mathbb{R}^\ell$ is as follows: 
\begin{eqnarray*}
F_{\alpha}\circ f=
\biggl(F_{1}\circ f+\sum_{j=1}^{m}\alpha_{1j}f_j, \ldots,
F_{\ell}\circ f+\sum_{j=1}^{m}\alpha_{\ell j}f_j\biggr),
\end{eqnarray*}
where $f=(f_1,\ldots ,f_m)$. 
Since there is the natural identification 
$\mathcal{L}(\mathbb{R}^{m},\mathbb{R}^{\ell})=(\mathbb{R}^{m})^\ell$, 
in order to prove \cref{@f}, 
it is sufficient to show that 
there exists a subset $\Sigma$ 
with Lebesgue measure zero of $(\mathbb{R}^m)^{\ell}$ 
such that for any $\alpha \in (\mathbb{R}^m)^{\ell}-\Sigma$, 
the mapping ${}_{s}j^{r}(F_{\alpha}\circ f):
N^{(s)}\to {}_{s}J^{r}(N,\mathbb{R}^\ell)$ is transverse to the given 
 modular submanifold $W$. 

Let $H_{\Lambda}:\mathbb{R}^\ell \to \mathbb{R}^\ell $
be the linear isomorphism defined by 
\begin{eqnarray*}
H_{\Lambda}(X_1,\ldots ,X_{\ell})=(X_1,\ldots ,X_{\ell})\Lambda,
\end{eqnarray*}
where $\Lambda=
(\lambda_{ij})_{1\le i\le \ell, 1\le j\le \ell}$ is an $\ell \times \ell$ regular matrix. 
The composition of $H_{\Lambda}$ and $F_{\alpha}\circ f$ is as follows:
{\footnotesize
\begin{eqnarray*}
H_{\Lambda}\circ F_{\alpha}\circ f
=\biggl(
\sum_{k=1}^{\ell}\biggl(F_{k}\circ f+\sum_{j=1}^{m}\alpha_{kj}f_j\biggr)\lambda_{k1},
\ldots,
\sum_{k=1}^{\ell}\biggl(F_{k}\circ f+\sum_{j=1}^{m}\alpha_{kj}f_j\biggr)\lambda_{k\ell}
\biggr)
\\
=\biggl(
\sum_{k=1}^{\ell}\biggl(F_{k}\circ f\biggr)\lambda_{k1}
+\sum_{j=1}^{m}\biggl(\sum_{k=1}^{\ell}\lambda_{k1}\alpha_{kj}\biggr)f_j,
\ldots,
\sum_{k=1}^{\ell}\biggl(F_{k}\circ f\biggr)\lambda_{k \ell}
+\sum_{j=1}^{m}\biggl(\sum_{k=1}^{\ell}\lambda_{k \ell}\alpha_{kj}\biggr)f_j
\biggr).
\end{eqnarray*}
}

Set $GL(\ell)=\{B\mid B:\ell \times \ell$ matrix, ${\rm det} B\not=0\}$.
Let $\varphi:GL(\ell)\times (\mathbb{R}^m)^{\ell}\to 
GL(\ell)\times (\mathbb{R}^m)^{\ell}$ 
be the mapping as follows:
{\footnotesize
\begin{eqnarray*}
{ } & { } &\varphi(\lambda_{11},\lambda_{12},\ldots ,\lambda_{\ell \ell},
\alpha_{11},\alpha_{12},\ldots ,\alpha_{\ell m})
\\&&=
\biggl(\lambda_{11},\lambda_{12},\ldots ,\lambda_{\ell \ell},
\sum_{k=1}^{\ell}\lambda_{k1}\alpha_{k1}, 
\sum_{k=1}^{\ell}\lambda_{k2}\alpha_{k1},\ldots ,
\sum_{k=1}^{\ell}\lambda_{k\ell}\alpha_{k1},
\\&&
\sum_{k=1}^{\ell}\lambda_{k1}\alpha_{k2}, 
\sum_{k=1}^{\ell}\lambda_{k2}\alpha_{k2},\ldots ,
\sum_{k=1}^{\ell}\lambda_{k\ell}\alpha_{k2}
,\ldots ,
\sum_{k=1}^{\ell}\lambda_{k1}\alpha_{km}, 
\sum_{k=1}^{\ell}\lambda_{k2}\alpha_{km},\ldots ,
\sum_{k=1}^{\ell}\lambda_{k\ell}\alpha_{km}
\biggr).
\end{eqnarray*}
}For the proof of \cref{@f}, 
it is the key to show that $\varphi$ is a $C^\infty$ diffeomorphism. 
In order to show that $\varphi$ is a $C^\infty$ diffeomorphism, 
for any point $(\Lambda' ,\alpha')\in GL(\ell)\times (\mathbb{R}^m)^{\ell}$ 
of the target space of $\varphi$, 
we will find $(\Lambda ,\alpha)$ satisfying 
$\varphi(\Lambda ,\alpha) =(\Lambda' ,\alpha')$, 
where $\Lambda=(\lambda_{11},\lambda_{12},\ldots ,\lambda_{\ell \ell})$, 
$\Lambda'=(\lambda'_{11},\lambda'_{12},\ldots ,\lambda'_{\ell \ell})$, 
$\alpha=(\alpha_{11},\alpha_{12},\ldots , \alpha_{\ell m})$, 
and $\alpha'=(\alpha'_{11},\alpha'_{12},\ldots , \alpha'_{m\ell})$. 
Hence, it is sufficient to find $(\Lambda ,\alpha)$ satisfying 
\begin{eqnarray*}
\lambda_{ij}&=&\lambda'_{ij} 
~(1\leq i \leq \ell, 1\leq j \leq \ell), 
\\
\sum_{k=1}^{\ell}\lambda_{ki}\alpha_{kj}&=&\alpha'_{ji}  
~(1\leq i \leq \ell, 1\leq j \leq m).
\end{eqnarray*}
Therefore, for any $j$ $(1\leq j \leq m)$, 
we get 
\begin{eqnarray*}
\sum_{k=1}^{\ell}\lambda'_{k1}\alpha_{kj}=\alpha'_{j1}, ~
\sum_{k=1}^{\ell}\lambda'_{k2}\alpha_{kj}=\alpha'_{j2},
\ldots , ~
\sum_{k=1}^{\ell}\lambda'_{k\ell}\alpha_{kj}=\alpha'_{j\ell}. 
\end{eqnarray*}
Thus, for any $j$ $(1\leq j \leq m)$, we have the following:
\begin{eqnarray*}
\left(
\begin{array}{cccccc}
\lambda'_{11} & \cdots & \lambda'_{\ell 1} \\ 
\vdots & \ddots & \vdots \\ 
\lambda'_{1\ell} & \cdots & \lambda'_{\ell \ell} \\ 
\end{array}
\right)
\left(
\begin{array}{cccccc}
\alpha_{1j}\\
\vdots \\
\alpha_{\ell j}
\end{array}
\right)
=
\left(
\begin{array}{cccccc}
\alpha'_{j1}\\
\vdots \\
\alpha'_{j \ell}
\end{array}
\right).
\end{eqnarray*}
Since the matrix 
\begin{eqnarray*}
\left(
\begin{array}{cccccc}
\lambda'_{11} & \cdots & \lambda'_{\ell 1} \\ 
\vdots & \ddots & \vdots \\ 
\lambda'_{1\ell} & \cdots & \lambda'_{\ell \ell} \\ 
\end{array}
\right)
\end{eqnarray*}
is regular, for any $j$ $(1\leq j \leq m)$, $\alpha_{1j},\ldots ,\alpha_{\ell j}$ 
can be expressed by rational functions of 
$\lambda'_{11}, \ldots , \lambda'_{\ell \ell}$, 
$\alpha'_{j1},\ldots ,\alpha'_{j \ell}$. 
Therefore, there exists the inverse mapping $\varphi^{-1}$ 
and we see that $\varphi^{-1}$ is of class $C^{\infty}$. 
Hence, the mapping $\varphi$ is a $C^{\infty}$ diffeomorphism. 

Next, let  $\widetilde{f}:U \to \mathbb{R}^{m+\ell}$ be the mapping 
as follows:
\begin{eqnarray*}
\widetilde{f}(x_1,\ldots,x_m)=(F_1(x_1,\ldots ,x_m),\ldots ,
F_{\ell}(x_1,\ldots ,x_m),x_1,\ldots ,x_m).
\end{eqnarray*}
It is clearly seen that $\widetilde{f}$ is an embedding. 
Since $f:N\to U$ is an embedding, 
the mapping $\widetilde{f}\circ f:N \to \mathbb{R}^{m+\ell}$ 
is also an embedding:
\begin{eqnarray*}
\widetilde{f}\circ f=(F_1\circ f,\ldots ,F_{\ell}\circ f,f_1,\ldots ,f_m).
\end{eqnarray*}
In order to prove \cref{@f}, the following lemma is important. 
The following lemma is the special case of \cref{@m}. 
\begin{lemma}[\cite{GP}] \label{@gp}
Let $N$ be a manifold of dimension $n$. 
Let $\widetilde{f}\circ f:N\to \mathbb{R}^{m+\ell}$ be an embedding. 
If $W$ is a modular submanifold of ${}_{s}J^{r}(N,\mathbb{R}^\ell)$, 
then there exists a subset $\Sigma$ with Lebesgue measure zero 
of $\mathcal{L}(\mathbb{R}^{m+\ell},\mathbb{R}^{\ell})$ 
such that for any $\Pi 
\in \mathcal{L}(\mathbb{R}^{m+\ell},\mathbb{R}^{\ell})-\Sigma$, 
the mapping 
${}_{s}j^{r}(\Pi \circ \widetilde{f}\circ f) : 
N^{(s)}\to {}_{s}J^{r}(N,\mathbb{R}^\ell)$ 
is transverse to $W$. 
\end{lemma}

From \cref{@gp}, 
there exists a subset $\Sigma$ with Lebesgue measure zero of 
$\mathcal{L}(\mathbb{R}^{m+\ell},\mathbb{R}^{\ell})$
such that for any $\Pi
\in \mathcal{L}(\mathbb{R}^{m+\ell},\mathbb{R}^{\ell})-\Sigma$, 
the mapping ${}_{s}j^{r}(\Pi \circ (\widetilde{f}\circ f)) : 
N^{(s)}\to {}_{s}J^{r}(N,\mathbb{R}^\ell)$ 
is transverse to $W$. 

There is the natural identification 
$\mathcal{L}(\mathbb{R}^{m+\ell},\mathbb{R}^{\ell})=\mathbb{R}^{\ell(m+\ell)}$. 
Thus, we identify the target space $GL(\ell)\times (\mathbb{R}^m)^{\ell}$ of the mapping $\varphi$ with 
an open submanifold of $\mathcal{L}(\mathbb{R}^{m+\ell},\mathbb{R}^{\ell})$. 
Since the intersection  $(GL(\ell)\times (\mathbb{R}^m)^{\ell})\cap \Sigma$ is 
a subset with Lebesgue measure zero 
of $GL(\ell)\times (\mathbb{R}^m)^{\ell}$ and 
the mapping $\varphi^{-1}$ is of class $C^\infty$, 
it follows that 
$\varphi^{-1}((GL(\ell)\times (\mathbb{R}^m)^{\ell})\cap \Sigma)$ is a subset with Lebesgue measure zero of 
$GL(\ell)\times (\mathbb{R}^m)^{\ell}$. 
For any $(\Lambda, \alpha) \in 
GL(\ell)\times (\mathbb{R}^m)^{\ell}$, 
let $\Pi_{(\Lambda, \alpha)}:
\mathbb{R}^{m+\ell} \to \mathbb{R}^{\ell} $ 
be the linear mapping defined by  $\varphi(\Lambda, \alpha)$ 
as follows:
\begin{eqnarray*}
{ } & { } &\Pi 
_{(\Lambda, \alpha)}(X_1,\ldots ,X_{m+\ell})\\
{ } & = &(X_1,\ldots ,X_{m+\ell})
\left(
\begin{array}{cccccc}
\lambda_{11} & \cdots & \lambda_{1\ell} \\ 
\vdots & \ddots & \vdots \\ 
\lambda_{\ell 1} & \cdots & \lambda_{\ell \ell} \\ 
\displaystyle{
\sum_{k=1}^{\ell}\lambda_{k1}\alpha_{k1}} 
& \cdots & 
\displaystyle{
\sum_{k=1}^{\ell}\lambda_{k\ell}\alpha_{k1}} \\ 
\vdots & \ddots & \vdots \\ 
\displaystyle{
\sum_{k=1}^{\ell}\lambda_{k1}\alpha_{km}} 
& \cdots & 
\displaystyle{
\sum_{k=1}^{\ell}\lambda_{k\ell}\alpha_{km}} \\ 
\end{array}
\right).
\end{eqnarray*}
Then, we have the following: 
{\footnotesize
\begin{eqnarray*}
{ } & { } &\Pi_
{(\Lambda, \alpha)}\circ \widetilde{f}\circ f
\\
{ } & = &\biggl(
\sum_{k=1}^{\ell}\biggl(F_{k}\circ f\biggr)\lambda_{k1}
+\sum_{j=1}^{m}\biggl(\sum_{k=1}^{\ell}\lambda_{k1}\alpha_{kj}\biggr)f_j,
\ldots,
\sum_{k=1}^{\ell}\biggl(F_{k}\circ f\biggr)\lambda_{k\ell}
+\sum_{j=1}^{m}\biggl(\sum_{k=1}^{\ell}\lambda_{k\ell}\alpha_{kj}\biggr)f_j
\biggr)
\\
{ } & = &H_{\Lambda}\circ F_{\alpha}\circ f.
\end{eqnarray*}
}Therefore, for any $(\Lambda, \alpha) \in 
GL(\ell)\times (\mathbb{R}^m)^{\ell}
- \varphi^{-1}((GL(\ell)\times (\mathbb{R}^m)^{\ell})\cap \Sigma)$, 
it follows that 
${}_sj^r(\Pi_{(\Lambda, \alpha)}\circ \widetilde{f}\circ f)$ 
$(={}_sj^r(H_\Lambda \circ F_{\alpha}\circ f))$ is 
transverse to $W$. 
Since the mapping $H_{\Lambda}$ is a diffeomorphism, 
we see that ${}_s j^r(F_{\alpha}\circ f)$ is transverse to $W$. 

Let $\widetilde{\Sigma}$ be a subset 
consisting of $\alpha \in (\mathbb{R}^m)^{\ell}$ 
such that ${}_s j^r(F_\alpha \circ f)$ is not transverse to $W$. 
In order to prove \cref{@f}, 
it is sufficient to show that $\widetilde{\Sigma}$ is a subset 
with Lebesgue measure zero of $(\mathbb{R}^m)^\ell$. 
Suppose that $\widetilde{\Sigma}$ is not a subset 
with Lebesgue measure zero of $(\mathbb{R}^m)^\ell$. 
Then, $GL(\ell)\times \widetilde{\Sigma}$ is not a subset
with Lebesgue measure zero 
 of $GL(\ell)\times (\mathbb{R}^m)^{\ell}$. 
For any $(\Lambda, \alpha) \in GL(\ell)\times \widetilde{\Sigma}$, 
since ${}_sj^r(F_{\alpha}\circ f)$ is not transverse to $W$ and 
the mapping $H_{\Lambda}$ is a diffeomorphism, 
${}_sj^r(H_{\Lambda}\circ F_{\alpha}\circ f)$ 
is not transverse to $W$. 
This contradicts to the claim that 
$\varphi^{-1}((GL(\ell)\times (\mathbb{R}^m)^{\ell})\cap \Sigma)$ is a subset 
with Lebesgue measure zero of $GL(\ell)\times (\mathbb{R}^m)^{\ell}$. 
\hfill $\Box$
\section{Applications of the main result}\label{section 4}
\subsection{Introduction of generalized distance-squared mappings}
In this subsection, the definition of generalized distance-squared mappings 
and the motivation to investigate the mappings are given. 
Moreover, for the sake of reader's convenience, the main properties 
of generalized distance-squared mappings are also reviewed 
(for more details on properties of generalized distance-squared mappings, refer to \cite{D}, \cite{L}, \cite{G2}, \cite{G1}). 

Let $i$, $j$ be positive integers, and 
let $p_i=(p_{i1}, p_{i2}, \ldots, p_{im})$  $(1\le i\le \ell)$ 
(resp., $A=(a_{ij})_{1\le i\le \ell, 1\le j\le m}$) 
be points of $\mathbb{R}^m$ 
(resp., an $\ell\times m$ matrix with non-zero entries).
Set $p=(p_1,p_2,\ldots,p_{\ell})\in (\mathbb{R}^m)^{\ell}$. 
Let $G_{(p, A)}:\mathbb{R}^m \to \mathbb{R}^\ell$ be the mapping 
defined by 
{\small 
\[
G_{(p, A)}(x)=\left(
\sum_{j=1}^m a_{1j}(x_j-p_{1j})^2, 
\sum_{j=1}^m a_{2j}(x_j-p_{2j})^2, 
\ldots, 
\sum_{j=1}^m a_{\ell j}(x_j-p_{\ell j})^2
\right), 
\] }where $x=(x_1, x_2, \ldots, x_m)\in \mathbb{R}^m$. 
The mapping $G_{(p, A)}$ is called a {\it generalized distance-squared mapping}, 
and the $\ell$-tuple of points $p=(p_1,p_2,\ldots ,p_{\ell})\in (\mathbb{R}^m)^{\ell}$ 
is called the {\it central point} 
of the generalized distance-squared mapping $G_{(p,A)}$. 

A {\it distance-squared mapping} $D_p$ 
(resp., {\it Lorentzian distance-squared mapping}
$L_p$) is the mapping $G_{(p,A)}$ 
satisfying that each entry of $A$ is $1$ 
(resp., $a_{i1}=-1$ and $a_{ij}=1$ $(j\ne 1)$). 

In \cite{D} (resp., \cite{L}), 
a classification result on distance-squared mappings $D_p$ 
(resp., Lorentzian distance-squared mappings
$L_p$) is given. 

In \cite{G1}, 
a classification result on 
generalized distance-squared mappings of the plane into the plane 
is given. 
If the rank of $A$ is equal to two,  
then a generalized distance-squared mapping 
having a generic central point 
is a stable mapping of which any singular point is 
a fold point 
except one cusp point 
(for details on fold points and cusp points, refer to \cite{plane to plane}).
If the rank of $A$ is equal to one, 
then a generalized distance-squared mapping 
having a generic central point is $\mathcal{A}$-equivalent 
to the normal form of definite fold mapping $(x_1,x_2)\to (x_1,x_2^2)$. 
Since the normal form of definite fold mapping is proper, 
it is easily shown that the mapping is stable 
by Mather's characterization theorem of stable proper mappings 
given in \cite{mather5}. 

In \cite{G2}, 
a classification result on 
generalized distance-squared mappings 
of $\mathbb{R}^{m+1}$ into $\mathbb{R}^{2m+1}$ 
is given. 
If the rank of $A$ is equal to $m+1$, 
then a generalized distance-squared mapping 
having a generic central point 
is $\mathcal{A}$-equivalent to 
the mapping called the normal form of Whitney umbrella as follows:
\[(x_1,\ldots ,x_{m+1})\mapsto (x_1^2,x_1x_2,\ldots ,x_1x_{m+1},x_2,\ldots ,x_{m+1}).\]
The normal form of Whitney umbrella is proper and stable. 
If the rank of $A$ is less than $m+1$, 
then a generalized distance-squared mapping 
having a generic central point 
is $\mathcal{A}$-equivalent to the inclusion as follows:
\[(x_1,\ldots ,x_{m+1})\mapsto (x_1,\ldots ,x_{m+1},0,\ldots ,0).\]
The inclusion is proper and stable. 

Quadratic polynomial mappings have been investigated for particular pairs 
of dimensions $(m,\ell)$ in different areas of mathematics, 
and there exists a vast literature on the subject. 
For example, in \cite{quadratic1} (resp., \cite{quadratic2}), 
a classification of quadratic polynomial mappings 
of the plane into the plane 
(resp., of the plane into the $n$ dimensional space) is given. 
In \cite{net}, nets of quadrics are investigated. 

Hence, as a research on quadratic polynomial mappings, 
it is natural to investigate the stability of 
generalized distance-squared mappings on submanifolds and 
the stability of generalized distance-squared mappings 
of $\mathbb{R}^m \to \mathbb{R}^\ell$ 
in the cases that $(m, \ell)$ is neither 
$(2, 2)$ nor $(k+1,2k+1)$, where $k$ is 
a positive integer. 

We have another original motivation. 
Height functions and distance-squared functions have been investigated 
in detail so far, and they are useful tools 
in the applications of singularity theory to differential geometry 
(for example, see \cite{CS} and \cite{Izumiyabook}). 
A mapping in which each component is a height function is nothing but a projection. 
Projections as well as height functions or distance-squared functions have been 
investigated so far. 
In \cite{GP}, the stability of projections on submanifolds is investigated.

On the other hand, 
a mapping in which each component 
is a distance-squared function is a distance-squared mapping. 
Moreover, the notion of a generalized distance-squared mapping is 
an extension of that of a distance-squared mapping. 
Therefore, we investigate the generalized distance-squared mappings 
as well as projections 
on submanifolds from the view point of the stability 
(see \cref{@Glis} and \cref{remark G}). 
\subsection{Applications of \cref{@f} to $G_{(p,A)} : \mathbb{R}^m\to \mathbb{R}^\ell$ }
As an application of \cref{@f}, we have the following. 
\begin{proposition}\label{@G}
Let $N$ be a manifold of dimension $n$. 
Let $f:N\to \mathbb{R}^m$ be an embedding. 
Let $A=(a_{ij})_{1\leq i \leq \ell, 1\leq j \leq m}$ 
be an $\ell \times m$ matrix with non-zero entries. 
If $W$ is a modular submanifold of ${}_{s}J^{r}(N,\mathbb{R}^\ell)$,
then there exists a subset $\Sigma$ with Lebesgue measure zero 
of $(\mathbb{R}^m)^\ell$ 
such that 
for any $p=(p_1,\ldots,p_\ell) \in  (\mathbb{R}^m)^\ell-\Sigma$, 
$G_{(p,A)}\circ f:N\to \mathbb{R}^\ell$ 
is transverse with respect to $W$. 
\end{proposition}
{\it Proof.}\qquad 
Let $H:\mathbb{R}^{\ell} \to \mathbb{R}^{\ell}$ 
be the diffeomorphism of the target for deleting constant terms. 
The composition $H\circ G_{(p,A)}:\mathbb{R}^{m}\to \mathbb{R}^{\ell}$ 
is given as follows:
\begin{eqnarray*}
H\circ G_{(p, A)}(x)=\left(
\sum_{j=1}^m a_{1j}x_j^2-2\sum_{j=1}^m a_{1j}p_{1j}x_j, 
\ldots, 
\sum_{j=1}^m a_{\ell j}x_j^2-2\sum_{j=1}^m a_{\ell j}p_{\ell j}x_j
\right), 
\end{eqnarray*}
where $x=(x_1,\ldots ,x_m)$. 

Let $\psi :(\mathbb{R}^m)^\ell \to \mathcal{L}(\mathbb{R}^{m},\mathbb{R}^{\ell})$ 
be the mapping defined by 
\begin{eqnarray*}
\psi (p_{11},p_{12},\ldots ,p_{\ell m})=-2(a_{11}p_{11},a_{12}p_{12},\ldots ,
a_{\ell m}p_{\ell m}). 
\end{eqnarray*}
Note that there exists the natural identification 
$\mathcal{L}(\mathbb{R}^{m},\mathbb{R}^{\ell})=(\mathbb{R}^m)^\ell$. 
Since $a_{ij}\not =0$ for any $i$, $j$ $(1\leq i \leq \ell$, $1\leq j \leq m)$, 
it is clearly seen that $\psi$ is a $C^\infty$ diffeomorphism. 

Set $F_{i}(x)=\sum_{j=1}^m a_{ij}x_j^2$ $(1\leq i \leq \ell)$ and $F=(F_1,\ldots ,F_\ell)$.  
From \cref{@f}, 
there exists a subset $\Sigma$ with Lebesgue measure zero 
of $\mathcal{L}(\mathbb{R}^{m},\mathbb{R}^{\ell})$ 
such that 
for any $\pi \in \mathcal{L}(\mathbb{R}^{m},\mathbb{R}^{\ell})-\Sigma$, 
the mapping ${}_{s}j^{r}(F_{\pi}\circ f):N^{(s)}\to {}_{s}J^{r}(N,\mathbb{R}^\ell)$ 
is transverse to $W$. 
Since $\psi^{-1}:\mathcal{L}(\mathbb{R}^{m},\mathbb{R}^{\ell}) \to 
(\mathbb{R}^m)^\ell$ is a $C^\infty$ mapping, 
$\psi^{-1}(\Sigma)$ is a subset 
with Lebesgue measure zero of $(\mathbb{R}^m)^\ell$. 
For any $p\in (\mathbb{R}^m)^\ell-\psi^{-1}(\Sigma)$, 
we have $\psi(p)\in \mathcal{L}(\mathbb{R}^{m},\mathbb{R}^{\ell})-\Sigma$. 
Hence, for any $p\in (\mathbb{R}^m)^\ell-\psi^{-1}(\Sigma)$, 
the mapping 
${}_{s}j^{r}(H\circ G_{(p, A)}\circ f):N^{(s)}\to {}_{s}J^{r}(N,\mathbb{R}^\ell)$ 
is transverse to $W$. 
Then, since $H:\mathbb{R}^{\ell} \to \mathbb{R}^{\ell}$ 
is a diffeomorphism, 
${}_{s}j^{r}(G_{(p, A)}\circ f):N^{(s)}\to {}_{s}J^{r}(N,\mathbb{R}^\ell)$ 
is transverse to $W$. 
\hfill 
$\Box$ 

\medskip 
As in \cref{section 2}, 
from \cref{@G}, we get two applications. 
\begin{corollary}
\label{TB}
Let $N$ be a manifold of dimension $n$. 
Let $f:N\to \mathbb{R}^m$ be an embedding. 
Let $A=(a_{ij})_{1\leq i \leq \ell, 1\leq j \leq m}$ 
be an $\ell \times m$ matrix with non-zero entries. 
Then, there exists a subset $\Sigma$ with Lebesgue measure zero 
of $(\mathbb{R}^m)^\ell$ 
such that 
for any $p=(p_1,\ldots,p_\ell) \in  (\mathbb{R}^m)^\ell-\Sigma$, 
$G_{(p,A)}\circ f:N\to \mathbb{R}^\ell$ 
is transverse with respect to the Thom-Boardman 
varieties.
\end{corollary} 
\begin{corollary}
\label{@Glis}
Let $N$ be a manifold of dimension $n$. 
Let $f:N\to \mathbb{R}^m$ be an embedding. 
Let $A=(a_{ij})_{1\leq i \leq \ell, 1\leq j \leq m}$ 
be an $\ell \times m$ matrix with non-zero entries. 
If a dimension pair $(n,\ell)$ is in the nice dimensions, 
then there exists a subset $\Sigma$ with Lebesgue measure zero 
of $(\mathbb{R}^m)^\ell$ 
such that 
for any $p=(p_1,\ldots,p_\ell) \in  (\mathbb{R}^m)^\ell-\Sigma$, 
the composition $G_{(p,A)}\circ f:N\to \mathbb{R}^\ell$ 
is locally infinitesimally stable. 
\end{corollary} 
\begin{remark}\label{remark G}
{\rm 
\begin{enumerate}
\item 
Suppose that the mapping 
$G_{(p,A)}\circ f:N\to \mathbb{R}^{\ell}$ is proper 
in \cref{@Glis}. 
Then, the local infinitesimal stability of $G_{(p,A)}\circ f$ 
implies the stability of it (see \cite{mather5}). 
\item 
Suppose that $N=\mathbb{R}^m$ and $f:\mathbb{R}^m \to \mathbb{R}^m$ is 
the identify. 
From \cref{@Glis}, 
it is clearly seen that if $(m,\ell)$ is in the nice dimensions, 
then there exists a subset $\Sigma$ with Lebesgue measure zero of $(\mathbb{R}^m)^\ell$
such that 
for any $p=(p_1,\ldots,p_\ell)\in (\mathbb{R}^m)^\ell-\Sigma$, 
the mapping $G_{(p,A)}:\mathbb{R}^m\to \mathbb{R}^\ell$ 
is locally infinitesimally stable. 
This is an application of (I5) in \cref{section 1}.  
\end{enumerate}
}
\end{remark}

\section*{Acknowledgements}
The author is most grateful to the anonymous reviewer for carefully reading 
the first manuscript of this paper and for giving invaluable suggestions. 
The author is grateful to Takashi Nishimura and 
Stanis\l aw Janeczko
for their kind advice. 
The author is supported by JSPS KAKENHI Grant Number 16J06911. 


\begin{thebibliography}{99}
\bibitem{Arnold}V.~I.~Arnol'd,  
\textit{Singularities of smooth mappings}, 
(Russian) 
Uspehi Mat. Nauk, \textbf{23} (1968), 3--44. 

\bibitem{Boardman}J.~M.~Boardman,  
\textit{Singularities of differentiable maps}, 
Inst. Hautes \'Etudes Sci. Publ. Math., \textbf{33} (1967), 21--57.





\bibitem{CS}J.~W.~Bruce and P.~J.~Giblin, 
\textit{Curves and singularities $($second edition$)$}, 
Cambridge University Press, Cambridge, 1992. 


\bibitem{brucekirk}J.~W.~Bruce and  N.~P.~Kirk, 
\textit{Generic projections of stable mappings}, 
Bull. London Math.
Soc., \textbf{32} (2000), 718--728.

\bibitem{net}S.~A.~Edwards and C.~T.~C.~Wall, 
\textit{Nets of quadrics and deformations of $\Sigma^{3(3)}$ singularities}, 
Math. Proc. Cambridge Philos. Soc., 
\textbf{105} (1989), no.~1, 109--115. 

\bibitem{quadratic1}M.~Farnik and Z.~Jelonek, 
\textit{On quadratic polynomial mappings of the plane}, 
Linear Algebra Appl., 
\textbf{529} (2017), 441--456. 

\bibitem{quadratic2}M.~Farnik, Z.~Jelonek and P.~Migus,  
\textit{On quadratic polynomial mappings from the plane into the $n$ 
dimensional space}, available from ResearchGate. 






\bibitem{GG}M.~Golubitsky and V.~Guillemin, 
\textit{Stable mappings and their singularities}, 
Graduate Texts in Mathematics \textbf{14}, Springer, New York, 1973.


\bibitem{D}S.~Ichiki and T.~Nishimura, 
\textit{Distance-squared mappings}, 
Topology Appl., 
\textbf{160} (2013), 1005--1016.   

\bibitem{L}S.~Ichiki and T.~Nishimura, 
\textit{Recognizable classification of Lorentzian distance-squared 
mappings}, J.~Geom.~Phys., \textbf{81} (2014), 62--71.   

\bibitem{G2}S.~Ichiki and T.~Nishimura, 
\textit{Generalized distance-squared mappings of 
$\mathbb{R}^{n+1}$ into $\mathbb{R}^{2n+1}$}, 
Contemp. Math., Amer. Math. Soc., Providence, RI, \textbf{675} (2016),
 121--132.

\bibitem{G1}S.~Ichiki, T.~Nishimura, R.~Oset Sinha and M.~A.~S.~Ruas, 
\textit{Generalized distance-squared mappings of the plane into the plane}, 
Adv. Geom., \textbf{16} (2016), 189--198.


\bibitem{Izumiyabook}
S.~Izumiya, M.~C.~Romero~Fuster, M.~A.~S.~Ruas and F.~Tari,  
\textit{Differential geometry from a singularity theory viewpoint}, 
World Scientific Publishing Co. Pte. Ltd., Hackensack, NJ, 2016.





\bibitem{mather5}J.~N.~Mather, 
\textit{Stability of $C^\infty$ mappings V. Transversality}, 
{Adv.~in~Math.}, 
\textbf{4} (1970), 301--336. 

\bibitem{mather6}J.~N.~Mather, 
\textit{Stability of $C^\infty$ mappings V\!I. The nice dimensions}, 
{Lecture
Notes in Math.}, 
\textbf{192} (1971), 207--253. 



\bibitem{TB}J. N. Mather, 
\textit{On Thom-Boardman singularities}, 
Dynamical systems (Proc. Sympos., Univ. Bahia, Salvador, 1971), 
Academic Press, New York, 1973, 233--248. 

\bibitem{GP}J. N. Mather, \textit{Generic projections}, 
Ann. of Math., (2) \textbf{98} (1973), 226--245.


\bibitem{semi}A.~A.~du~Plessis and  C.~T.~C~Wall, 
\textit{Generic projections in the semi-nice dimensions}, 
Compositio Math. 
\textbf{135} (2003), 179--209.

\bibitem{Thom}R.~Thom, 
\textit{Les singularit\'es des applications diff\'erentiables}, 
Ann. Inst. Fourier, Grenoble, \textbf{6} (1955--1956), 43--87.


\bibitem{plane to plane}H.~Whitney, 
\textit{On singularities of mappings of euclidean spaces. I. 
Mappings of the plane into the plane}, 
Ann. of Math., (2), 
\textbf{62} (1955), 374--410.   
  
\end{thebibliography}
\end{document}